\documentclass{amsart}

\usepackage{amsmath}
\usepackage{amscd}
\usepackage{amssymb}

\input xy
\xyoption{all}

\newcommand{\pd}{\partial}

\newcommand{\bC}{{\mathbb C}}

\newcommand{\bP}{{\mathbb P}}

\newcommand{\bR}{{\mathbb R}}
\newcommand{\bZ}{{\mathbb Z}}

\newcommand{\cM}{{\mathcal M}}
\newcommand{\cO}{{\mathcal O}}
 \newcommand{\cR}{{\mathcal R}}

\DeclareMathOperator{\Hom}{Hom}

\newtheorem{Example}{Example}

\newtheorem{theorem/definition}{Theorem/Definition}[section]

\newtheorem{Conjecture}{Conjecture}

\theoremstyle{remark}

\theoremstyle{definition}

\newcommand{\be}{\begin{equation}}
\newcommand{\ee}{\end{equation}}
\newcommand{\bea}{\begin{eqnarray}}
\newcommand{\eea}{\end{eqnarray}}
\newcommand{\ben}{\begin{eqnarray*}}
\newcommand{\een}{\end{eqnarray*}}
\newcommand{\bet}{\begin{equation}
\begin{split}}
\newcommand{\eet}{\end{split}
\end{equation}}

\begin{document}

\title
{Crepant Resolutions, Quivers and GW/NCDT Duality}
\author{Jian Zhou}
\address{Department of Mathematical Sciences\\Tsinghua University\\Beijing, 100084, China}
\email{jzhou@math.tsinghua.edu.cn}

\begin{abstract}
We propose a conjecture that relates some local Gromov-Witten invariants
of some crepant resolutions of Calabi-Yau $3$-folds with isolated singularities
with some Donaldson-Thomas type invariants of the moduli spaces of representations
of some quivers with potentials.
\end{abstract}
\maketitle

\section{Introduction}

Singular Calabi-Yau $3$-folds and their crepant resolutions are
very interesting both in algebraic geometry and in string theory.
Not only do they provide numerous interesting examples of resolutions
of singularities,
but also the study of their invariants is very rich.
Given a crepant resolution $\pi: Y \to X$ of an affine Calabi-Yau $3$-fold
with an isolated singularity,
the local Gromov-Witten invariants and the derived category of coherent sheaves of $Y$
are of particular interest.
They can be studied both from the algebrogeometric point of view
and the string theoretical point of view,
and the interactions are very crucial for some of the recent progresses.
Gromov-Witten invariants are defined mathematically as
intersection numbers on the moduli spaces of stable maps to $Y$,
they correspond to some correlators in type IIA closed string theory;
objects in derived category of coherent sheaves are complexes of coherent sheaves up to
some equivalence relations,
they correspond to type IIB D-branes.
Let us recall how D-branes lead us to quivers with superpotentials functions.
Physically Type B D-brane at a singular point can decay into a collection of stable fractional D-branes.
Each constituent D-brane may appear with a multiplicity,
and there are open strings between them.
In this way one can associates a quiver gauge theory to it.
Mathematically,
this is described by a noncommutative resolution \cite{VdB}:
one can often find a quiver with superpotential $(Q, W)$ such that there is an equivalence:
$$D^b(Y) \cong D^b(Q, W),$$
where $D^b(Y)$ is the derived category of coherent sheaves on $Y$,
and $D^b(Q, W)$ is the derived category of quiver representations of $Q$ constrained
by the relations given by the superpotential function.
The path algebra of the quiver $Q$ with relations given by $W$
is noncommutative,
called a Calabi-Yau algebra \cite{Boc, Gin}.
One can define Donaldson-Thomas type invariants of the moduli spaces
of semistable representations of $(Q, W)$.
We conjecture that they are related to the Gromov-Witten invariants of $Y$.
For more precise formulation,
See Conjecture 1 in Section 3.

The rest of the paper is organized as follows.
In Section 2 we will recall some examples of crepant resolutions of affine Calabi-Yau
$3$-folds with isolated singularities
and the quivers with superpotentials associated to them.
In Section 3 we will present our conjecture.

\section{Crepant Resolutions of Calabi-Yau $3$-Folds and Quivers with Potentials}

In this section we will recall some well-known examples of crepant resolutions and their associated
quivers with superpotentials.

\subsection{The case of $\bP^1$ as exceptional sets}

By a result due to Laufer \cite{Lau},
when the exceptional set is $\bP^1$,
its normal bundle is isomorphic to one of the following three bundles:
$\cO(-1) \oplus \cO(-1)$, $\cO(-2) \oplus \cO$, and $\cO(-3)\oplus \cO(1)$.
The first case can be realized by the well-known resolved conifold,
and the their local Gromov-Witten invariants are well-known \cite{Gop-Vaf, Fab-Pan};
the other two case are realized by Laufer's examples \cite{Lau},
for their local Gromov-Witten invariants, and more generally,
that of $\cO(k) \oplus \cO(-k-2) \to \bP^1$ for $k \geq -1$
see \cite{Yan-Zho}.

The resolved conifold can be obtained by gluing two copies of $\bC^3$,
with linear coordinates $(x, y_1, y_2)$ and $(w, z_1, z_2)$ respectively,
by the following formula for change of coordinates:
\begin{align*}
z_1 & = xy_1, &
z_2 & = x y_2, &
w & = \frac{1}{x}.
\end{align*}
Let
\begin{align*}
v_1 & = z_1 = xy_1, & v_2 & = z_2 = xy_2, \\
v_3 & = wz_1 = y_1, & v_4 & = wz_2 = y_2.
\end{align*}
They satisfy:
\be
v_1v_4 -  v_2 v_3 = 0.
\ee
Hence we get a contraction
$\pi: Y \to X :=\{(v_1, v_2, v_3, v_4) \in \bC^4: \;\; v_1v_4 - v_2v_3 = 0\}$.
The associated quiver is (cf. \cite{Asp-Kat} and the references therein):
$$
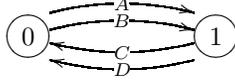
\xymatrix{
*++[o][F-]{0}  \ar@/^/[rr]|-{B} \ar@/^/@<1.5ex>[rr]|-{A}   & &
*++[o][F-]{1}  \ar@/^/[ll]|-{C}  \ar@/^/@<1.5ex>[ll]|-{D} }
$$
with superpotential function:
$$W = BCAD-ACBD.$$
This case has been studied by Nagao and Nakajima \cite{Nag-Nak}.
Note the conifold and the resolved conifold are toric Calabi-Yau $3$-folds.
In general we will consider toric affine Calabi-Yau$3$-folds
that correspond to the polygons with integral vertices in $\bZ^2$.
A toric crepant resolution corresponds to a subdivision into triangles with integral vertices
and of area $1/2$.
If two subdivisions differ only by changing the diagonal of a parallelogram,
then the corresponding crepant resolutions differ by a flop.
For example,
the conifold corresponds to the square with vertices
$(0, 0)$, $(1, 0)$, $(0, 1)$ and $(1, 1)$.
It has two subdivisions:
$$
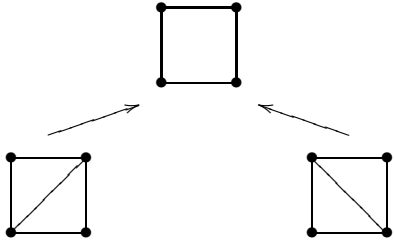
\xy
(0,0); (10,0), **@{-}; (10, 10), **@{-}; (0, 10), **@{-};(0, 0), **@{-};
(0,0)*+{\bullet}; (10,0)*+{\bullet};(10, 10)*+{\bullet};(0, 10)*+{\bullet};
(-20,-20); (-10,-20), **@{-}; (-10, -10), **@{-}; (-20, -10), **@{-}; (-20, -20), **@{-};
(-20, -20); (-10,-10), **@{-};
(-20,-20)*+{\bullet}; (-10,-20)*+{\bullet};(-10, -10)*+{\bullet};(-20, -10)*+{\bullet};
(20,-20); (30,-20), **@{-}; (30, -10), **@{-}; (20, -10), **@{-};(20, -20), **@{-};
(20, -10); (30, -20), **@{-};
(20,-20)*+{\bullet}; (30,-20)*+{\bullet};(30, -10)*+{\bullet};(20, -10)*+{\bullet};
(-15,-7); (-3,-3), **\dir{-} ?>* \dir{>}; (25,-7); (13,-3), **\dir{-} ?>* \dir{>};
\endxy
$$

Now we come to Laufer's examples.
The first class of examples is the space $Y$ obtained by gluing two copies of $\bC^3$,
with linear coordinates $(x, y_1, y_2)$ and $(w, z_1, z_2)$ respectively,
by the following rules for change of coordinates:
\begin{align*}
z_1 & = x^2y_1 + xy_2^k, &
z_2 & = y_2, &
w & = \frac{1}{x}.
\end{align*}
Define a map from $Y$ to $\bC^4$ by:
\begin{align*}
v_1 & = z_2 = y_2, &
v_2 & = z_1 = x^2y_1+xy_2^k, \\
v_3 & = wz_1 = xy_1 + y_2^k, &
v_4 & = w^2 - wz_2^k = y_1.
\end{align*}
The image is the affine variety $X$ defined by
$v_2v_4-v_3^2+v_3v_1^k=0$, or by a change of coordinates:
\be
u_1^2+ u_2^2 + u_2^2 + u_4^{2k} = 0.
\ee
The exceptional set is given by $(x, y_1=0, y_2)$ and $(w, z_1 = 0, z_2=0)$.
Its normal bundle is isomorphic to $\cO(-1) \oplus \cO(-1) \to \bP^1$ when $k=1$,
and to $\cO\oplus \cO(-2) \to \bP^1$ when $k > 1$.
For $k > 1$,
the associated quiver is  (cf. \cite{Asp-Kat} and the references therein):
$$
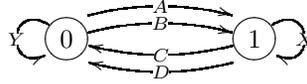
\xymatrix{
*++[o][F-]{0} \ar@(ul,dl)[]|{Y} \ar@/^/[rr]|-{B} \ar@/^/@<1.5ex>[rr]|-{A}   & &
*++[o][F-]{1} \ar@(ur,dr)[]|{X} \ar@/^/[ll]|-{C}  \ar@/^/@<1.5ex>[ll]|-{D} }
$$
with superpotential function
$$W = - (-1)^{n(n-1)/2}X^{n+1}-(-1)^{n(n-1)/2}Y^{n+1}-XAC+XBD-YCA+YDB.$$

These examples are not toric,
nevertheless there is a natural $2$-torus action defined as follows:
\ben
&& (t_1, t_2) \cdot (x, y_1, y_2) \mapsto (t_1^{-1}t_2^k x, t_1y_1, t_2y_2), \\
&& (t_1, t_2) \cdot (w, z_1, z_2) \mapsto (t_1 t_2^{-k} w, t_1^{-1} t_2^{2k} z_1, t_2 z_2),
\een
where $t_1, t_2 \in \bC^*$.
Using this action one can define and compute the local Grmov-Witten invariants
as in \cite{Yan-Zho}.

The second class of examples is the space $Y$ obtained by gluing two copies of $\bC^3$,
with linear coordinates $(x, y_1, y_2)$ and $(w, z_1, z_2)$ respectively,
by the following formula for change of coordinates:
\begin{align*}
z_1 & = x^3y_1 + y_2^2 + x^2y_2^{2n+1}, &
z_2 & = \frac{y_2}{x}, &
w & = \frac{1}{x}.
\end{align*}
Define a map $Y \to \bC^4$ by:
\begin{align*}
v_1 & = z_1 && = x^3y_1+xy_2^{2n+1}, \\
v_2 & = w^2z_1 - z_2^2 && = xy_1+y_2^{2n+1}, \\
v_3 & = w^3z_1 -wz_2^2-z_1^nz_2 && = y_1 + \frac{1}{x}[y_2^{2n+1}-y_2(x^3y_1+y_2^2+x^2y_2^{2n+1})^n], \\
v_4 & = w^2z_1z_2-z_2^3 - wz_2^{n+1} && = y_1y_2 + \frac{1}{x} [y_2^{2n+2}-(x^3y_1+y_2^2+x^2y_2^{2n+1})^{n+1}].
\end{align*}
The image is the affine variety defined by:
\be
v_4^2+ v_2^3 - v_1v_3-  v_1^{2n+1}v_2 = 0.
\ee
The associated quiver does not seem to be known in the literature.
There is a natural $\bC^*$-action on $Y$ defined as follows:
\begin{eqnarray*}
&& t \cdot (x, y_1, y_2) = (t^{1-2n}x, t^{6n+1}y_1, t^2y_2), \\
&& t \cdot (w, z_1, z_2) = (t^{2n-1} w, t^4z_1, t^{2n+1} z_2),
\een
where $t\in \bC^*$.
Using this action one can define and compute the local Gromov-Witten invariants
as in \cite{Yan-Zho}.

\subsection{The case of a string of rational curves as exceptional sets}

The case of crepant resolutions of a toric Calabi-Yau $3$-fold
whose exceptional set is a string of rational curves
has been studied by Nagao \cite{Nag}.
There are two cases for the corresponding lattice polygons:
Case 1. The triangle with vertices $(0,0)$, $(2,0)$ and $(0, 2)$:
$$
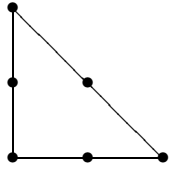
\xy
(0,0); (20,0), **@{-}; (0, 20), **@{-}; (0, 0), **@{-};
(0,0)*+{\bullet}; (10,0)*+{\bullet};(10, 10)*+{\bullet}; (0, 10)*+{\bullet};
(20,0)*+{\bullet};(0, 20)*+{\bullet};
\endxy
$$
It has four subdivisions related by flops:
$$
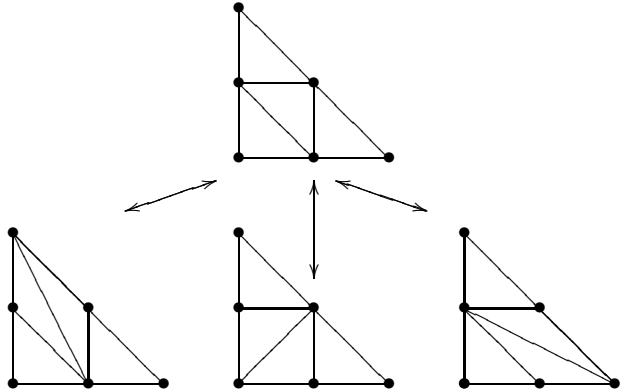
\xy
(0,0); (20,0), **@{-}; (0, 20), **@{-}; (0, 0), **@{-};
(0,10);(10,10), **@{-}; (10,0), **@{-}; (0,10), **@{-};
(0,0)*+{\bullet}; (10,0)*+{\bullet};(10, 10)*+{\bullet}; (0, 10)*+{\bullet};
(20,0)*+{\bullet};(0, 20)*+{\bullet};
(0,-30); (20,-30), **@{-}; (0, -10), **@{-}; (0, -30), **@{-};
(0,-20); (10,-20), **@{-}; (10, -30), **@{-}; (0,-30), **@{-}; (10,-20), **@{-};
(0,-30)*+{\bullet}; (10,-30)*+{\bullet};(10, -20)*+{\bullet}; (0, -20)*+{\bullet};
(20,-30)*+{\bullet};(0, -10)*+{\bullet};
(-30,-30); (-10,-30), **@{-}; (-30, -10), **@{-}; (-30, -30), **@{-};
(-30,-20); (-20, -30), **@{-}; (-30,-10),**@{-}; (-20,-20), **@{-}; (-20,-30), **@{-};
(-30,-30)*+{\bullet}; (-20,-30)*+{\bullet};(-20, -20)*+{\bullet}; (-30, -20)*+{\bullet};
(-10,-30)*+{\bullet};(-30, -10)*+{\bullet};
(30,-30); (50,-30), **@{-}; (30, -10), **@{-}; (30, -30), **@{-};
(30,-20); (40,-20), **@{-}; (50,-30), **@{-}; (30,-20), **@{-}; (40,-30), **@{-};
(30,-20)*+{\bullet}; (50,-30)*+{\bullet};(30, -10)*+{\bullet}; (30, -30)*+{\bullet};
(30,-30)*+{\bullet}; (40, -30)*+{\bullet};  (40, -20)*+{\bullet};
(-15,-7); (-3,-3), **\dir{-} ?>* \dir{>}; (-3,-3); (-15,-7), **\dir{-} ?>* \dir{>};
(25,-7); (13,-3), **\dir{-} ?>* \dir{>}; (13,-3); (25,-7), **\dir{-} ?>* \dir{>};
(10, -3); (10, -16), **\dir{-} ?>* \dir{>}; (10, -16); (10, -3), **\dir{-} ?>* \dir{>};
\endxy
$$
Case 2. The trapezoid with vertices $(0, 0)$, $(N_0, 0)$, $(N_1, 1)$ and $(0,1)$,
where $N_0 \geq N_1 > 0$ are positive integers.
For example:
$$
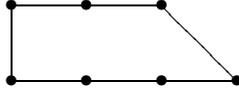
\xy
(0,0); (30,0), **@{-}; (20, 10), **@{-}; (0,10), **@{-}; (0, 0), **@{-};
(0,0)*+{\bullet}; (10,0)*+{\bullet}; (20,0)*+{\bullet};(30, 0)*+{\bullet};
(20,10)*+{\bullet}; (10, 10)*+{\bullet}; (0, 10)*+{\bullet};
\endxy
$$
Again there are many different ways to subdivide it into some triangles of area $1/2$.

In both cases,
one can use localization to define and compute their local Gromov-Witten invariants,
and express the results in terms of the topological vertex \cite{AKMV, LLLZ}.
For the second case, see \cite{Iqb-Kas}.

\subsection{The case of surfaces as exceptional sets}

When the exceptional set of $\pi: Y \to X$ is an algebraic surface,
it can only be a del Pezzo surface \cite{Rei},
i.e.
$\bP^2$ or $\bP^1 \times \bP^1$ blown up at $k < 9$ points.
For the associated quivers and their superpotentials,
see e.g. \cite{Asp-Fid, Asp-Mel, Fen-Fra-Han-He, He} and the references therein.
For example,
for $Y = \cO_{\bP^2}(-3)$,
the associated quiver is
$$
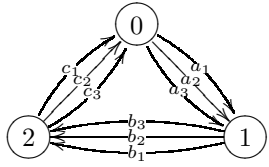
\xymatrix{
 & *++[o][F-]{0} \ar[rd]|-{a_2} \ar@/^/[rd]|-{a_1}  \ar@/_/[rd]|-{a_3} & \\
*++[o][F-]{2}  \ar[ru]|-{c_2} \ar@/^/[ru]|-{c_1}  \ar@/_/[ru]|-{c_3} & &
*++[o][F-]{1}  \ar[ll]|-{b_2} \ar@/^/[ll]|-{b_1}  \ar@/_/[ll]|-{b_3} }
$$
with superpotential function
$$
W = a_1b_2c_3-a_1b_3c_2+a_2b_3c_1-a_2b_1c_3+a_3b_1c_2-a_3b_2c_1.
$$
By choosing different strongly exceptional collections,
it is possible to arrive at different quivers with superpotentials.
For example,
in the del Pezzo $3$ case,
there are four different possibilities (see e.g. \cite{Fen-Fra-Han-He}).

There are many toric examples (see e.g. \cite{CKYZ}, Fig. 1):
$$ 
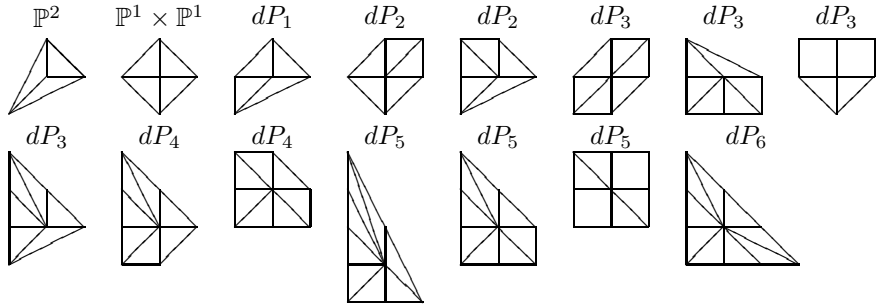
\xy
(0,0);(5,0), **@{-}; (0,5), **@{-}; (0,0), **@{-}; (-5,-5), **@{-}; (0,5),**@{-};(5,0),**@{-};(-5,-5),**@{-};
(0,8)*+{\bP^2};
(10,0); (15,5),**@{-}; (20,0),**@{-}; (15,-5), **@{-}; (10,0), **@{-}; (20,0),**@{-};
(15,5); (15,-5),**@{-}; (15,8)*+{\bP^1\times \bP^1};
(25,0); (30,5),**@{-}; (35,0),**@{-};(25,-5),**@{-}; (25,0), **@{-}; (35,0),**@{-};
(30,5); (30,0),**@{-}; (25,-5), **@{-};
(30,8)*+{dP_1};
(40,0); (45,5),**@{-}; (50,5), **@{-}; (50,0), **@{-}; (45,-5), **@{-}; (40,0), **@{-};
(40,0); (50,0), **@{-}; (45,5); (45, -5), **@{-}; (45,0); (50,5), **@{-};
(45,8)*+{dP_2};
(55,5); (60,5),**@{-}; (65,0), **@{-}; (55,-5), **@{-}; (55, 5), **@{-};
(55, 0); (65, 0),**@{-}; (55,5); (60, 0), **@{-};(60,5),**@{-}; (60,0); (55,-5)**@{-};
(60,8)*+{dP_2};
(70,-5); (70,0),**@{-}; (75,5),**@{-}; (80, 5),**@{-}; (80,0),**@{-}; (75,-5),**@{-};(70,-5),**@{-};
(80,5),**@{-}; (70,0); (80,0),**@{-}; (75,-5); (75,5),**@{-};
(75,8)*+{dP_3};
(85,-5); (85,5),**@{-}; (95,0),**@{-};(95,-5),**@{-};(85,-5),**@{-};(90,0),**@{-};(90,-5),**@{-};
(85,0); (95,0),**@{-}; (85,5); (95,-5),**@{-};
(90,8)*+{dP_3};
(100,5);(110,5), **@{-};(110,0),**@{-}; (105,-5), **@{-}; (100,0),**@{-}; (100,5),**@{-};
(100,0); (110,0),**@{-}; (105,5); (105,-5),**@{-};
(105,8)*+{dP_3};
(0,-8)*+{dP_3};
(-5,-10); (-5,-25),**@{-}; (5,-20),**@{-}; (-5,-10), **@{-}; (0,-20), **@{-}; (-5,-25), **@{-};
(-5,-15); (0,-20),**@{-}; (0,-15),**@{-}; (-5,-20); (5,-20),**@{-};
(10,-10);(10,-25),**@{-}; (15,-25),**@{-}; (20,-20),**@{-};(10,-10),**@{-}; (15,-20),**@{-};(10,-25),**@{-};
(10,-15);(15,-20),**@{-}; (10,-20); (20,-20),**@{-}; (15,-25); (15,-15),**@{-};
(15,-8)*+{dP_4};
(25,-10); (25,-20),**@{-}; (35, -20),**@{-}; (35,-15),**@{-}; (30,-10),**@{-}; (25, -10),**@{-}; (35,-20),**@{-};
(25,-15); (35,-15),**@{-}; (30,-10); (30, -20),**@{-}; (25,-20); (30, -15),**@{-};
(30,-8)*+{dP_4};
(40,-10); (40,-30), **@{-}; (50,-30),**@{-}; (40,-10), **@{-}; (45, -25),**@{-}; (40,-15),**@{-};
(40,-30),**@{-}; (45, -25), **@{-};
(40,-20); (50,-30), **@{-}; (45,-20); (45, -30), **@{-}; (40,-25); (45,-25),**@{-}; (45,-8)*+{dP_5};
(55,-10); (55,-25), **@{-}; (65,-25), **@{-}; (65,-20), **@{-}; (55,-10), **@{-};
(60,-20), **@{-}; (55, -25), **@{-};
(55,-15); (65,-25), **@{-}; (55,-20); (65,-20), **@{-}; (60,-15); (60,-25), **@{-};
(60,-8)*+{dP_5};
(70,-10);(70,-20), **@{-}; (80,-20), **@{-}; (80,-10), **@{-}; (70,-10), **@{-}; (80,-20), **@{-};
(70,-15);(80,-15), **@{-}; (75,-10); (75,-20), **@{-}; (75,-8)*+{dP_5};
(85,-10); (85,-25), **@{-}; (100,-25), **@{-};(85,-10), **@{-}; (90,-20), **@{-};(85,-25), **@{-};
(85,-15); (95,-25), **@{-}; (85,-20); (95,-20), **@{-}; (90,-15); (90,-25), **@{-};
(90,-20); (100,-25), **@{-}; (93,-8)*+{dP_6};
\endxy$$
where $dP_n$ means $\bP^2$ blown up at $n$ points.
The local Gromov-Witten invariants of these toric examples have been computed in \cite{Zho}.
For some nontoric cases,
see \cite{Dia-Flo, Kon-Min}.

\subsection{Crepant resolutions of $3$-dimensional Gorenstein singularities}

Recall there is a classification of finite subsgroups $\Gamma \subset SL(3, \bC)$ (cf. \cite{Mil}).
Let $X = \bC^3/\Gamma$.
By a case by case analysis,
one can obtain a crepant resolution $\widehat{\bC^3/\Gamma}$
of $\bC^3/\Gamma$ (see Roan \cite{Roa} and the references given there).
Many examples are toric Calabi-Yau $3$-fold,
so one can use localization to define and compute their local Gromov-Witten invariants,
and express the results in terms of the topological vertex \cite{AKMV, LLLZ}.

\begin{Example}
Let $\Gamma$ be generated by
$$\begin{pmatrix}
e^{2\pi i/3} & 0 & 0 \\ 0 & e^{2\pi i/3} & 0 \\ 0 & 0 & e^{2\pi i/3} \end{pmatrix}.$$
Then one can take $\widehat{\bC^3/\Gamma} \cong \kappa_{\bP^2}$.
It can be covered by three local coordinate patches $(u_i, v_i, w_i) \in \bC^3$,
$i=1, 2, 3$.
Let $(z_1, z_2, z_3)$ be linear coordinates on $\bC^3$.
We take
\begin{align*}
u_1 & = \frac{z_2}{z_1}, & v_1 & = \frac{z_3}{z_1}, & w_1 & = z_1^3, \\
u_2 & = \frac{z_1}{z_2}, & v_2 & = \frac{z_3}{z_2}, & w_2 & = z_2^3, \\
u_3 & = \frac{z_1}{z_3}, & v_3 & = \frac{z_2}{z_3}, & w_3 & = z_3^3.
\end{align*}
From these we get:
\begin{align*}
u_2 & = \frac{1}{u_1}, & v_2 & = \frac{v_1}{u_1}, & w_2 & = u_1^3w_1,
\end{align*}
and so on.
Let $(\bC^*)^3$ act on $\bC^3$ by
$$(t_1, t_2, t_3) \cdot (z_1, z_, z_3)
= (t_1z_1, t_2z_2, t_3z_3).$$
This induces the following action on the $\kappa_{\bP^2}$ given in the above local coordinate patches
by
\begin{align*}
u_1 & \mapsto \frac{t_2}{t_1} \cdot u_1, & v_1 & \mapsto \frac{t_3}{t_1} \cdot v_1, &
w_1 & \mapsto t_1^3 \cdot w_1, \\
u_2 & \mapsto \frac{t_1}{t_2} \cdot u_2, & v_2 & \mapsto \frac{t_3}{t_2} \cdot v_2, &
w_2 & \mapsto t_2^3 \cdot w_2, \\
u_3 & \mapsto \frac{t_1}{t_3} \cdot u_3, & v_3 & \mapsto \frac{t_2}{t_3} \cdot v_3, &
w_3 & \mapsto t_3^3 \cdot w_3.
\end{align*}
It is clear that that $\{(u_1, v_1, w_1) \in \bC^3\;|\; u_1v_1w_1 \neq 0\}$
is a dense open orbit of the torus action.
Each coordinate patch has a unique fixed point $(u_i, v_i, w_i) = (0, 0, 0)$.
Denote them by $p_i$.
The weight decomposition at these points are
$(-\alpha+\alpha_2) \oplus (-\alpha_1+\alpha_3) \oplus (3\alpha_1)$,
$(\alpha_1-\alpha_2) \oplus (-\alpha_2+\alpha_3) \oplus (3\alpha_2)$
$(\alpha_1-\alpha_3) \oplus (\alpha_2-\alpha_3) \oplus (3\alpha_3)$
respectively.
The information about the fixed points and the weight decompositions at the fixed
points of the torus action can be encoded in the GKM graph \cite{GKM}:
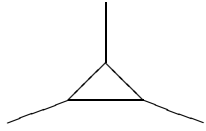
$$
\xy
(0,8); (0,0), **@{-}; (-5, -5), **@{-}; (5, -5), **@{-};  (0,0), **@{-};
(-5, -5); (-13,-8), **@{-}; (5, -5); (13, -8), **@{-};
\endxy$$
The toric information can also be encoded in the following lattice polygon:
$$
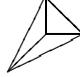
\xy
(0,0);(5,0), **@{-}; (0,5), **@{-}; (0,0), **@{-}; (-5,-5), **@{-}; (0,5),**@{-};(5,0),**@{-};(-5,-5),**@{-};
\endxy $$
This coincides with one of the cases we have seen in the previous subsection.
\end{Example}

\begin{Example}
In general,
for an odd positive integer $2n+1$,
let $\Gamma$ be generated by the matrix
$$\begin{pmatrix}
e^{2\pi i/(2n+1)} & 0 & 0 \\ 0 & e^{2\pi i/(2n+1)} & 0 \\ 0 & 0 & e^{-4\pi i/(2n+1)}
\end{pmatrix},
$$
one can explicitly write down a crepant resolution using $2n+1$ coordinate patches.
For example, when $n=3$,
one can take local coordinates:
\begin{align*}
&& (\frac{z_1}{z_3^3}, \frac{z_2}{z_3^3}, z_3^7) &&  \\
(\frac{z_2}{z_1}, \frac{z_3^3}{z_1}, \frac{z_1^3}{z_3^2}) &&   && (\frac{z_1}{z_2}, \frac{z_3^3}{z_2}, \frac{z_2^3}{z_3^2}) \\
(\frac{z_2}{z_1}, \frac{z_3^2}{z_1^3}, \frac{z_1^5}{z_3}) && && (\frac{z_1}{z_2}, \frac{z_3^2}{z_2^3}, \frac{z_2^5}{z_3}) \\
(\frac{z_2}{z_1}, \frac{z_3}{z_1^5}, z_1^7) && && (\frac{z_1}{z_2}, \frac{z_3}{z_2^5}, z_2^7).
\end{align*}
The corresponding GKM graph looks like:
$$
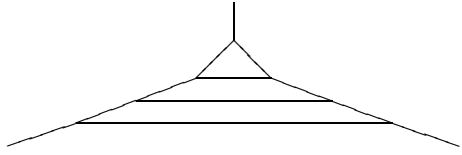
\xy
(0,5); (0,0), **@{-}; (-5, -5), **@{-}; (5, -5), **@{-};  (0,0), **@{-};
(-5, -5); (-13,-8), **@{-}; (13, -8), **@{-};(5, -5), **@{-};
(-13,-8); (-21, -11), **@{-}; (21, -11), **@{-}; (13, -8), **@{-};
(-21, -11); (-30, -14), **@{-}; (21, -11); (30, -14), **@{-};
\endxy$$
The toric information can also be encoded in the lattice polygon with vertices $(1,0)$, $(0,1)$,
$(-n,-n)$.
The crepant resolution corresponding to the subdivision given by adding vertices at
$(-k,-k)$, $k=0,1, \dots, n-1$.
See e.g. the case of $n=2$:
$$
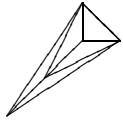
\xy
(0,5);(-5,-5), **@{-}; (5,0), **@{-}; (0,5), **@{-}; (-10,-10), **@{-}; (5,0), **@{-};
(0,5), **@{-}; (0,0), **@{-}; (5,0),**@{-};
(-10,-10); (0,0),**@{-};
\endxy $$
\end{Example}

\subsubsection{Three-dimensional McKay correspondence}
Another method to obtain the crepant resolution is via Nakamura's $\Gamma$-Hilbert schemes
\cite{Bri-Kin-Rei}.
This realizes a crepant resolution of $\bC^3/\Gamma$ as the moduli space
of $\Gamma$-cluster of points in $\bC^3$.
As a consequence, Bridgeland-King-Reid \cite{Bri-Kin-Rei}
proved by Fourier-Mukai transform an equivalence of derived categories:
\be
D^b(\widehat{\bC^3/\Gamma}) \cong D^b_\Gamma(\bC^3),
\ee
where $D^b(\widehat{\bC^3/\Gamma})$ denotes the derived category of bounded complexes
of coherent sheaves on $\widehat{\bC^3/\Gamma}$,
and $D^b_\Gamma(\bC^3)$ denotes the derived category of bounded complexes of $\Gamma$-equivariant
coherent sheaves on $\bC^3$.
This is a generalization of the McKay correspondence in two dimensions \cite{Gon-Spr, Kap-Vas}.
Define a noncommutative algebra
$$R_{\Gamma} = \bC[z_1,z_2,z_3] \ltimes \bC\Gamma,$$
where $\bC[z_1,z_2,z_3]$ is the algebra of polynomials on $\bC^3$,
$\bC \Gamma$ is the group ring of $\Gamma$,
and $G_{\Gamma}$ is their twisted product, in particular,
$$g \cdot p(z_1, z_2, z_3) = g(p(z_1,z_2,z_3)) \cdot g,$$
where $g \in \Gamma$, $p(z_1, z_2,z_3)$,
$g(p(z_1, z_2,z_3))$ is the action of $g$ on $p(z_1, z_2,z_3)$.
Let $R_\Gamma-mod$ be the category of finitely generated $R_\Gamma$-module and let
$D^b(R_{\Gamma}-mod)$ be its derived category.
Then we have
\be
D^b_\Gamma(\bC^3) \cong D^b(R_\Gamma).
\ee
The noncommutative algebra $R_\Gamma$ is an example of
noncommutative resolution \cite{VdB}.
Note the center of $R_\Gamma$ is $\bC[z_1, z_2, z_3]^{\Gamma}$,
the coordinate ring of $\bC^3/\Gamma$,
and so we have the natural inclusion $\bC[z_1,z_2,z_3]^\Gamma \hookrightarrow R_\Gamma$.

\subsubsection{Quivers with superpotentials}
Given an $R_\Gamma$-module $M$,
we can regard it first as a $\Gamma$-module and consider a decomposition:
$$M = \oplus_{\chi \in \hat{\Gamma}} M_\chi,$$
where $\hat{\Gamma}$ is the set of irreducible representations,
and $M_\chi$ is a direct sum of irreducible representations of characters $\chi$.
Then we consider $M$ as an $\bC[z_1, z_2, z_3]$-module and consider the maps of multiplications
by $z_1$, $z_2$ and $z_3$.
For simplicity we first assume $G$ is abelian,
and so all irreducible representations are $1$-dimensional.
Furthermore, assume $\bC z_i$ is an irreducible representation $\chi_{z_i}$ of $\Gamma$.
Then multiplication by $z_i$ induces a map
$$z_i: M_\chi \to  M_{\chi\otimes \chi_{z_i}}.$$
This inspires one to consider the following quiver (directed graph):
Assign a vertex to each $\chi$,
and assign an arrow from $\chi$ to $\chi \otimes \chi_{z_i}$ for each $i=1, 2, 3$.
Because multiplication by $z_i$ commutes with multiplication by $z_j$ when $i \neq j$,
some relations should be imposed for the quiver,
as is clear from the following example.

\begin{Example}
Let $\Gamma = \bZ_3$ and let $\omega$ a generator.
There are three irreducible representations $\chi_0$, $\chi_1$ and $\chi_2$:
$\chi_k(\omega) = e^{2k \pi i/3}$.
Suppose that
\be
\omega \cdot z_k = e^{2\pi i/3} \cdot z_k \omega.
\ee
The decomposition of the $R_{\Gamma}$-module is described by the following diagram:
$$
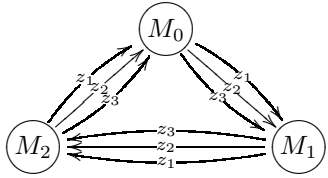
\xymatrix{
 & *++[o][F-]{M_0} \ar[rd]|-{z_2} \ar@/^/[rd]|-{z_1}  \ar@/_/[rd]|-{z_3} & \\
*++[o][F-]{M_2}  \ar[ru]|-{z_2} \ar@/^/[ru]|-{z_1}  \ar@/_/[ru]|-{z_3} & &
*++[o][F-]{M_1}  \ar[ll]|-{z_2} \ar@/^/[ll]|-{z_1}  \ar@/_/[ll]|-{z_3} }
$$
This leads to the quiver we have seen in the $\cO(-{\bP^2}(-3)$ case:
$$
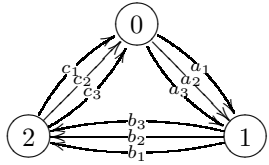
\xymatrix{
 & *++[o][F-]{0} \ar[rd]|-{a_2} \ar@/^/[rd]|-{a_1}  \ar@/_/[rd]|-{a_3} & \\
*++[o][F-]{2}  \ar[ru]|-{c_2} \ar@/^/[ru]|-{c_1}  \ar@/_/[ru]|-{c_3} & &
*++[o][F-]{1}  \ar[ll]|-{b_2} \ar@/^/[ll]|-{b_1}  \ar@/_/[ll]|-{b_3} }
$$
but now $a_1, a_2, a_3$ correspond to multiplication by $z_1$,
$b_1, b_2, b_3$  correspond to multiplication by $z_2$,
and $c_1, c_2, c_3$  correspond to multiplication by $z_3$.
Therefore,
we should impose the following relations:
\begin{align}
a_ib_j & = a_j b_i , & b_ic_j & = b_jc_i, & c_ia_j & = c_ja_i,
\end{align}
whenever $i \neq j$.
These $27$ relations can be encoded in the superpotential function:
$$
W = a_1b_2c_3-a_1b_3c_2+a_2b_3c_1-a_2b_1c_3+a_3b_1c_2-a_3b_2c_1.
$$
The relations are obtained by taking
$\frac{\pd}{\pd a_i}$, $\frac{\pd}{\pd b_i}$ and $\frac{\pd}{\pd c_i}$
for $i=1, 2,3$.
\end{Example}

For the construction of the quiver with potential $(Q^\Gamma, W_\Gamma)$
associated to a general $\Gamma \subset GL(3, \bC)$,
see Ginzburg \cite[\S 4.4]{Gin}.
Its path algebra is Morita equivalent to $R_\Gamma$.

\vspace{.2in}

\section{GW/NCDT Correspondence}

In last section,
we have seen many examples that to a crepant resolution $\pi: Y \to X$
of affine Calabi-Yau $3$-fold with isolated singularity,
one can associate a quiver $Q$ with superpotential function $W$,
such that
$$D^b(Y) \cong D^b(Q, W).$$
We have also seen that many examples of the crepant resolutions are toric,
and therefore one can use localization to define and compute their local Gromov-Witten invariants.
Each of these quivers have a special vertex $v_0$,
corresponding to $\cO_Y$ in the strong exceptional collection
that gives rise to $Q$.
Following Nagao-Nakajima \cite{Nag-Nak} and Nagao \cite{Nag},
we will consider the framed quiver $\tilde{Q}$ by adding a vertex $v_{\infty}$
and an arrow from $v_\infty$ to $v_0$, with the same superpotential function $W$.

\subsection{Moduli spaces of quiver representations}
Given a quiver $Q=(Q_0, Q_1)$,
where $Q_0$ denotes the set of vertices,
and $Q_1$ the set of arrows.
There are two maps $h, t: Q_1 \to Q_0$
which specify for each arrow its head and tail.

Let $\{W_v\}_{v \in Q_0}$ a collection of complex vector spaces,
one for each vertex of the quiver.
Its dimension vector of a representation is the vector $\alpha = (\dim W_v)_{v \in Q_0}$.
A representation of $Q$ with dimension vector $\alpha$ is a collection of
linear maps $\{\phi_a: W_{ta} \to W_{ha}\}_{a \in Q_1}$.
The dimension vector of a representation is the vector $\alpha = (\dim W_v)_{v \in Q_0}$.
They form a vector space
$$\cR(Q; \alpha)
: = \bigoplus_{a\in Q_1} \Hom(W_{ta}, W_{ha}),$$
called the representation space of $Q$ with dimension vector $\alpha$.
There is a natural action of the group
$$GL(\alpha):=\prod_{v \in Q_0} GL(W_v)$$
given by $(g \cdot \phi)_a = g_{ha}\phi_ag_{ta}^{-1}$.
To define the mdouli spaces of quiver representations,
one has to consider the stability conditions.
Given a collector of real numbers $\theta=(\theta_v)_{v \in Q_0}$,
it gives a stability condition if
\be
\sum_{v\in V_0} \theta_v \alpha_v = 0.
\ee
Given a stability condition $\theta$,
a $Q$-representation $(U, \phi)$ is $\theta$-semistable (stable)
if
\be
\sum_{v\in Q_0} \theta_v \dim U_v \geq (>) 0
\ee
for every subrepresentation.
Denote by $\cR_{\theta-ss}(Q;\alpha) \subset \cR(Q; \alpha)$ the set $\theta$-semistable
$Q$-representations of dimension vector $\alpha$,
and
$\cM_{\theta}(Q; \alpha): = \cR_{\theta-ss}(Q; \alpha)/GL(\alpha)$ the mdouli space.
If some relations coming from a potential function $\Phi$ are imposed on $Q$,
we denote $\cR(Q, \Phi;\alpha)\subset \cR(Q;\alpha)$ the subset of representations satisfying these relations,
and by $\cM_{\theta}(Q, \Phi; \alpha)$ the corresponding moduli space.
We conjecture that for suitable stability conditions and dimension vectors,
the moduli spaces give some crepant resolutions.

\subsection{Noncommutative Donaldson-Thomas invariants}
Now let $(Q, W)$ be the quiver with potential function associated with
a crepant resolution $\pi: Y \to X$ as in \S 2,
let $(\tilde{Q}, W)$ be the corresponding framed quiver.
We will consider dimension vectors $\tilde{\alpha}$ of $\tilde{Q}$ such that $\alpha_\infty = 1$.
For generic $\theta: Q_0 \to \bR$,
extend it to $\tilde{\theta}: Q_0\cup \{\infty\} \to \bR$
by $\tilde{\theta}_v= \theta_v$ for $v \in Q_0$ and
\be
\theta_{\infty} = - \sum_{v\in Q_0} \theta_v \alpha_v.
\ee
By the same argument as in Segal \cite{Seg} and Szendroi \cite{Sze},
$\cM_{\tilde{\theta}}(\tilde{Q}, W; \tilde{\alpha})$ has a symmetric perfect obstruction theory,
and one can use Behrend's constructible function \cite{Beh} to
define its Donaldson-Thomas type invariants.
We will refer to these invariants as noncommutative Donaldson-Thomas (NCDT) invariants.
Denote the corresponding partition function by $Z_{\tilde{\theta}}(\tilde{Q}, W)$.

As in \cite{Nag-Nak, Nag},
we consider all stability conditions.
We conjecture that they have some chamber structures,
related to the infinite root system associated to $Q$ by the Kac Theorem \cite{Kac}.
Recall one can associate a Cartan matrix to the quiver $Q$,
and use the Weyl reflections to define a root system.
There is an indecomposable representation of $Q^\Gamma$ of dimension vector $\alpha$
if and only if $\alpha$ is a positive root.
Furthermore,
if $\alpha$ is a positive real root,
then there is a unique representation of $Q$ of dimension vector $\alpha$;
if $\alpha$ is a positive imaginary root,
them such representations form a variety of dimension $(\alpha, \alpha)$.
The relations imposed by the potential function $W$ will eliminate  some positive roots
because the corresponding indecomposable representations do not satisfy the relations.
The remaining positive roots will then give rise to hypersurfaces that
divide the space of $\theta$'s into chambers.

\subsection{GW/NCDT correspondence conjecture}

Inspired by Crepant Resolution Conjecture \cite{Rua, Bry-Gra, Coa-Rua},
GW/DT correspondence conjecture \cite{MNOP},
and  recent work of Nagao and Nakajima \cite{Nag-Nak, Nag},
we make the following

\begin{Conjecture}
Suppose that $(Q, W)$ is a quiver with superpotential function associated to a crepant resolution
$\pi: Y \to X$,
where $X$ is an affine Calabi-Yau $3$-fold with an isolated singularity.
Let $(\tilde{Q}, W)$ be the framed quiver obtained from $(Q, W)$.
Then the generating series $Z^{GW}(Y)$ of local Gromov-Witten invariants of $Y$
can be identified with the generating series $Z^{NCDT}_{\theta}(\tilde{Q}, W)$ of
NCDT invariants of $(\tilde{Q}, W)$ for suitable
stability parameters $\theta$ in a chamber $C_{GW}$.
\end{Conjecture}

We will refer to this conjecture as the noncommutative crepant resolution conjecture
or GW/NCDT correspondence conjecture.
In joint work with Weiping Li in progress,
we are verifying our conjecture in the case of $\pi: \widehat{\bC^3/\bZ_3} \to \bC^3/\bZ_3$.

We conjecture that there is a chamber $C_0$ such that $Z_\theta(Q, W) = 1$
for $\theta \in C_0$,
and there is a sequence of wall-crossings to $C_{GW}$ such that each wall-crossing
changes the partition function by multiplying a factor so that $Z_{\theta}(\tilde{Q}, W)$
becomes an infinite product for $\theta \in C_{GW}$.
We expect this corresponds to the infinite product structure of
$Z^{GW}(Y)$ predicted by the Gopakumar-Vafa integrality \cite{Gop-Vaf}.
This may shed some lights on the Gopakumar-Vafa invariants by relating them to infinite root systems.

We speculate that our conjecture follows from more general results, such as the PT type invariants
associated to two different $t$-structures of a Calabi-Yau category
can be identified after suitable change of variables,
and the relationship between GW invariants, DT invariants and PT invariants of toric $Y$..

\vspace{.2in}

{\em Acknowledgements}.
The conjecture proposed in this paper started to take its shape during a visit to Hong Kong University
of Science and Technology,
where the author has greatly enjoyed the hospitality of the Department of Mathematics there.
The author thanks Professor Weiping Li for the invitation and for many helpful discussions.
He also thanks Yunfeng Jiang, Hiraku Nakajima, Yongbin Ruan, Jie Xiao and Bin Zhu for comments and discussions.
This research is partially supported by two NSFC grants (10425101 and 10631050)
and a 973 project grant NKBRPC (2006cB805905).


\begin{thebibliography}{99}

\bibitem{AKMV} M. Aganagic, A. Klemm, M. Mari$\tilde{n}$o, C. Vafa,
{\em The topological vertex},
Commun. Math. Phys. 254 (2005), no. 2, 425-478.

\bibitem{Asp-Fid}
P.S. ~Aspinwall, L. M.~Fidkowski
{\em Superpotentials for quiver gauge theories},
arXiv:hep-th/0506041.


\bibitem{Asp-Kat}
P.S.~Aspinwall, S.~Katz
{\em Computation of superpotentials for D-Branes},
arXiv:hep-th/0412209.

\bibitem{Asp-Mel}
P.S.~Aspinwall, I. V.~Melnikov,
{\em D-Branes on vanishing del Pezzo surfaces},
arXiv:hep-th/0405134.

\bibitem{Beh}
K.~Behrend,
{\em Donaldson-Thomas invariants via microlocal geometry},
arXiv:math/0507523.


\bibitem{Boc}
R.~Bocklandt,
{\em Graded Calabi Yau algebras of dimension 3},
J. Pure Appl. Algebra {\bf 212} (2008), no. 1, 14-32.

\bibitem{Bri-Kin-Rei}
T. Bridgeland, A. King, M. Reid,
{\em The McKay correspondence as an equivalence of derived categories},
J. Amer. Math. Soc. {\bf 14} (2001), 535-554.

\bibitem{Bry-Gra}
J.~Bryan, T.~Graber,
{\em The crepant resolution conjecture},
to appear in {\em Algebraic Geometry---Seatle 2005 Proceedings},
arXiv:math.AG/0610129.

\bibitem{CKYZ}
T.-M. Chiang, A. Klemm, S.-T. Yau, E. Zaslow,
{\em Local mirror symmetry: Calculations and interpretations},
Adv. Theor. Math. Phys. {\bf 3} (1999) 495-565, arXiv:hep-th/9903053.

\bibitem{Coa-Rua}
T.~Coates, Y.~Ruan,
{\em Quantum cohomology and crepant resolutions: A conjecture},
arXiv:math/0710.5901.

\bibitem{Dia-Flo}
D.-E.~Diaconescu, B.~Florea,
{\em The ruled vertex and nontoric del Pezzo surfaces},
arXiv:hep-th/0507240.

\bibitem{Fab-Pan}
C.~Faber, R.~Pandharipande,
{\em Hodge integrals and Gromov-Witten theory}.
Invent. Math. {\bf 139} (2000), no. 1, 173--199.


\bibitem{Fen-Fra-Han-He}
B. Feng, S. Franco, A. Hanany and Y. H. He,
{\em Unhiggsing the del Pezzo},
arXiv:hep-th/0209228.

\bibitem{Gin}
V. Ginzburg,
{\em Calabi-Yau algebras},
arXiv:math/0612139.

\bibitem{GKM}
M.~Goresky, R.~Kottwitz, R.~MacPherson,
{\em Equivariant cohomology, Koszul duality, and the localization theorem},
Invent. Math. {\bf 131} (1998), no. 1, 25-83.

\bibitem{Gon-Spr}
G. Gonzalez-Sprinberg, J.-L. Verdier,
{\em Construction g\'eom\'etrique de la correspondance de McKay},
Ann. Sci. \'Ecole Norm. Sup. {\bf 16} (1983) 409-449.

\bibitem{Gop-Vaf} R.Gopakumar,C.Vafa, {\em M-theory and Topological Strings-II}, arXiv:hep-th/9812127.

\bibitem{He}
Y.-H.~He,
{\em Lectures on D-branes, gauge theories and Calabi-Yau singularities},
arXiv:hep-th/0408142.

\bibitem{Iqb-Kas}
A.~Iqbal, A.-K.~ Kashani-Poor
{\em The vertex on a strip},
arXiv:hep-th/0410174.


\bibitem{Kac}
V. Kac,
{\em Root systems, representations of quivers and invariant theory}.
Invariant theory (Montecatini, 1982), 74-108,
Lecture Notes in Math., {\bf 996}, Springer, Berlin, 1983.

\bibitem{Kap-Vas}
M.~Kapranov, E. Vasserot,
{\em Kleinian singularities, derived categories and Hall algebras},
Math. Ann. {\bf 316} (2000), 565-576.

\bibitem{Kon-Min}
Y.~Konishi, S.~Minabe
{\em Local Gromov-Witten invariants of cubic surfaces via nef toric degeneration},
arXiv:math/0607187.

\bibitem{Lau}
H. Laufer,
{\em On $\bC\bP^1$ as an exceptional set}, Recent developments in several complex variables,
Ann. of Math. Stud. Princeton University Press, {\bf 100} (1981), 261-275.


\bibitem{LLLZ}
J. Li, C.-C. Liu, K. Liu, J. Zhou,
{\em A mathematical theory of the topological vertex},
arXiv:math/0408426.

\bibitem{MNOP} D. Maulik, N. Nekrasov, A. Okounkov, R. Pandharipande,
{\em Gromov-Witten theory and Donaldson-Thomas theory, I},
Comp. Math. {\bf 142} (2006), 1263-1285.

\bibitem{Mil}
G. A. Miller, H . F. Blichfeld, L. E. Dickson: Theory and applications of finite groups, John Wiley and
Sons. Inc. (1916).

\bibitem{Nag}
K.~ Nagao,
{\em Derived categories of small toric Calabi-Yau 3-folds and counting invariants},
arXiv:0809.2994.


\bibitem{Nag-Nak}
K.~Nagao, H.~Nakajima,
{\em Counting invariant of perverse coherent sheaves and its wall-crossing},
arXiv:0809.2992.

\bibitem{PT}
R. Pandharipande, R.P. Thomas,
{\em Curve counting via stable pairs
in the derived category}, arXiv:0707.2348.

\bibitem{Rei}
M.~Reid, {\em Canonical 3-folds}. In: Beauville, A. (ed.) G\'eom\'etrie alg6brique, Angers 1979,
pp. 273-310. The Netherlands: Sijthoff and Noordhoff 1990.

\bibitem{Roa}
S.-S.~Roan,
{\em  Minimal resolutions of Gorenstein orbifolds in dimension three},
Topology {\bf 35} (1996), no. 2, 489-508.


\bibitem{Rua}
Y.~Ruan,
{\em The cohomology ring of crepant resolutions of orbifolds},
in Gromov-Witten theory of spin curves and orbifolds, 117-126, Contemp. Math., 403.
Amer. Math. Soc., Providence, RI, 2006.

\bibitem{Seg}
E.~Segal,
{\em The A-infinity deformation theory of a point and the derived categories of local Calabi-Yaus},
arXiv:0702539.

\bibitem{Sze}
B.~Szendroi,
{\em Non-commutative Donaldson-Thomas theory and the conifold},
Geom.Topol. {\bf 12} (2008), 1171-1202,
arXiv:0705.3419.

\bibitem{VdB}
M. Van den Bergh,
{\em Non-commutative crepant resolutions}.
The legacy of Niels Henrik Abel, 749-770, Springer, Berlin, 2004.
arXiv:math.RA/0211064.

\bibitem{Yan-Zho}
F.~Yang, J.~Zhou,
{\em Local Gromov-Witten invariants and tautological sheaves on Hilbert schemes},
preprint, 2008.

\bibitem{Zho}
J.~Zhou,
{\em Localizations on moduli spaces and free field realizations of Feynman rules},
arXiv:math/0310283.

\end{thebibliography}
\end{document}